\documentclass[11pt]{article}
\usepackage[a4paper,margin=26mm]{geometry}
\usepackage{amsmath,amssymb,amsthm,mathtools,microtype}
\usepackage[colorlinks=true,linkcolor=blue,citecolor=blue,urlcolor=blue]{hyperref}
\newtheorem{theorem}{Theorem}
\newtheorem{lemma}{Lemma}
\newtheorem{proposition}{Proposition}
\newtheorem*{krtheorem}{Theorem (Kaiser--Rabinovich)}
\theoremstyle{remark}
\newcommand{\HD}{\operatorname{HD}}
\newcommand{\RR}{\mathbb R}
\newcommand{\F}{\mathcal F}
\newcommand{\C}{\mathcal C}
\newcommand{\D}{\mathcal D}
\newcommand{\R}{\mathcal R}
\title{On Exact Hadwiger-Debrunner Numbers}
\author{Chaya Keller\thanks{Institute for the Theory of Computing,
Faculty of Computer and Information Science,
Ben-Gurion University of the Negev,
Be'er Sheva 84105, Israel and
School of Computer Science, Ariel University, Israel. \texttt{chayak@bgu.ac.il}
}}
\date{}
\begin{document}
\maketitle
\vspace{-1.7em}
\begin{abstract}
Hadwiger and Debrunner proved in 1957 that a finite family of compact
convex sets in $\RR^d$ with the $(p,q)$-property can be pierced by
$p-q+1$ points whenever $q>(d-1)p/d+1$. The bound $p-q+1$ is optimal.
We prove that this optimal piercing bound is obtained already whenever
$q\ge p^{1/3+\varepsilon}$ in the plane or
$q\ge p^{\sqrt{1-2/d}+\varepsilon}$ in high dimensions.
This is the first
range in which exact Hadwiger-Debrunner numbers are obtained since the 1957 work of Hadwiger
and Debrunner.

We also show that the exponents in these ranges
cannot be improved without improving the currently known
upper-bound exponents for weak $\varepsilon$-nets.
\end{abstract}

\section{Introduction}

\subsection{Background}

For integers $p\ge q$, a family $\F$ (with $|\F|\geq p$) has the \emph{$(p,q)$-property}
if among every $p$
members there are $q$ with a common point. A \emph{transversal}, or
\emph{piercing set}, is a set of points meeting every member of
$\F$; its minimum size is denoted by $\tau(\F)$. 

Helly's theorem says that a $(d+1,d+1)$-family of convex sets in $\RR^d$ is pierced by one
point. Hadwiger and Debrunner \cite{HD57} asked whether, more
generally, the $(p,q)$-property for a family of compact convex sets in $\RR^d$, guarantees a transversal whose size is bounded in terms of
$p$, $q$, and $d$ alone
for every $p\ge q\ge d+1$. Alon and Kleitman proved this in their
celebrated $(p,q)$-theorem \cite{AK92,AK97}. 

We write
$\HD_d(p,q)$ for the smallest universal bound on the transversal
size, the \emph{Hadwiger--Debrunner number}.
The lower bound
\[
 \HD_d(p,q)\ge p-q+1
 \tag{1}\label{eq:standardlower}
\]
is immediate: take $p-q$ pairwise disjoint sets and $q$ further
sets with a common point, disjoint from the first sets. Hadwiger
and Debrunner proved that equality holds whenever
\[
 q>\frac{d-1}{d}p+1.
 \tag{2}\label{eq:classical}
\]
Extending the range in which the optimal value $p-q+1$ is attained
has remained a basic problem, even though $\HD_d(p,q)$ is known
to be finite for every $q \geq d+1$.

\paragraph{Weak $\varepsilon$-nets.}
The Alon--Kleitman proof uses the weak $\varepsilon$-net theorem
of Alon, B\'ar\'any, F\"uredi, and Kleitman \cite{ABFK92}.
For every $0<\varepsilon\le1$ and every $d\ge2$, there is a bound
$f_d(\varepsilon)$ such that every finite nonempty multiset $X$
of points in $\RR^d$ admits a set $\mathcal N$ of at most $f_d(\varepsilon)$
points meeting every convex set that contains at least
$\varepsilon|X|$ points of $X$. Such a set $\mathcal N$ is called a
\emph{weak $\varepsilon$-net} for $X$. We take $f_d(\varepsilon)$
to be the least such bound. 

The planar bound $O(\varepsilon^{-2})$ of \cite{ABFK92}, and the
higher-dimensional bound $\widetilde O_d(\varepsilon^{-d})$ of
\cite{CEGG+95}, stood for more than two decades. Rubin
\cite{RubinPlane,RubinHigher} proved that, for every $\sigma>0$,
\[
 f_d(\varepsilon)=O_{d,\sigma}(\varepsilon^{-\rho_d-\sigma}),
 \qquad
 \rho_d=
 \begin{cases}
 3/2,&d=2,\\
 2.558,&d=3,\\
 3.48,&d=4,\\
 (d+\sqrt{d^2-2d})/2,&d\ge5.
 \end{cases}
 \tag{3}\label{eq:rubin}
\]
In every dimension $d\ge2$, Bukh,
Matou\v{s}ek, and Nivasch \cite{BMN11} proved that
\[
f_d(\varepsilon)=\Omega_d\!\left(
 \varepsilon^{-1}\log^{d-1}(1/\varepsilon)\right).
 \tag{4}\label{eq:BMN}.
\]
For $d\ge3$, very recently Bhore, Browne, and Zheng \cite[Theorem~1]{BBZ26}
proved the stronger lower bound
$f_d(\varepsilon)=\Omega_d(\varepsilon^{-\gamma_d})$,
where $\gamma_d=\sqrt{d/2}+O(1)$.

\paragraph{The nearly optimal range.}
Keller, Smorodinsky, and Tardos \cite{KST18} substantially improved
the quantitative $(p,q)$ bounds. Their proof uses a combinatorial
dichotomy: for $p'<p$, $q'<q$, either $\F$ satisfies the
$(p',q')$-property, or there is $S\subseteq\F$ with $|S|=p'$
and no intersecting $q'$-tuple. In the latter case
$\F\setminus S$ has the $(p-p',q-q'+1)$-property. This allows
one to pass recursively to smaller parameters.

Substituting (\ref{eq:rubin}) into the weak-net estimates in
\cite{KST18} gives the nearly optimal bound
\[
 \HD_d(p,q)\le p-q+2
 \quad\text{when }q\ge p^{1-1/\rho_d+\eta}
 \text{ and }p\text{ is sufficiently large}.
 \tag{5}\label{eq:nearintro}
\]
Actually, the proof in \cite{KST18} gives the bound $p-q+2$ in the narrower
range $q\ge p^{(d-1)/d+\eta}$. Adapting that proof to
Rubin's improved weak $\varepsilon$-net bounds yields
(\ref{eq:nearintro}), as we show in Appendix~\ref{sec:nearsharp}.
It seems unlikely, however, that the proof technique of \cite{KST18}
can sharpen the bound from $p-q+2$ to the exact value $p-q+1$.

However, the exact bound $p-q+1$ has been obtained over wider ranges for
some special classes of convex sets. For example, 
\cite{KS20} proved it for axis-parallel boxes when
$q>C_d\log^{d-1}p$. This proof relies on the fact that
every pairwise intersecting family of axis-parallel boxes has
a common point. But even in the plane, a pairwise intersecting
family of convex sets can have empty total intersection, so
extending this approach to general convex sets requires a
different geometric argument.

\subsection{Our results}
In Theorem~\ref{thm:main}, we prove that the optimal
Hadwiger--Debrunner number $p-q+1$ is attained over a wide range
of $q$. The proof uses two main tools. In Section~2.1, we pass
to sequences of convex sets, which allows us to apply a theorem
of Kaiser and Rabinovich \cite{KR99}. In Section~2.2, we develop
a new dichotomy that combines this argument with the known bound
$p-q+2$.

\begin{theorem}\label{thm:main}
For every fixed $d\ge2$ and $\eta>0$, there is
$p_0=p_0(d,\eta)$ such that
\[
 p\ge p_0,\qquad p\ge q\ge d+1,\qquad
 q\ge p^{1-1/\rho_d+\eta}
 \quad\Longrightarrow\quad
 \boxed{\HD_d(p,q)=p-q+1.}
 \tag{6}\label{eq:main}
\]
In particular, the planar range is $q\ge p^{1/3+\eta}$.
For $d\ge5$, the exponent $1-1/\rho_d$ equals
$\sqrt{1-2/d}$.
\end{theorem}

Our second result relates the range in which the optimal value
$p-q+1$ is attained to quantitative bounds in the weak
$\varepsilon$-net theorem.
\begin{theorem}\label{thm:correspondence}
For every fixed $d\ge2$ and $s\ge1$, the following are equivalent.
\begin{enumerate}
\item For every $\sigma>0$,
$f_d(\varepsilon)=O_{d,s,\sigma}(\varepsilon^{-s-\sigma})$
as $\varepsilon\rightarrow0$.
\item For every $\eta>0$ and all sufficiently large $p$,
$\HD_d(p,q)=p-q+1$ for all 
$q\ge p^{1-1/s+\eta}$.
\end{enumerate}
\end{theorem}

Thus, any polynomial improvement
of the threshold in Theorem~\ref{thm:main} would improve the
currently known weak-net exponent, and conversely,
every improvement of a weak-net exponent improves the exact
$(p,q)$ range.

\paragraph{The widest range for the exact bound.}
How far can the range in which $\HD_d(p,q)=p-q+1$ be extended?
In the plane, a natural candidate is $q\ge C\log p$. Indeed, as observed in
\cite[Eq.~(1)]{KS21}, the lower bound construction for weak
$\varepsilon$-nets due to Bukh, Matou\v{s}ek, and Nivasch \cite{BMN11}
implies\footnote{Set $\varepsilon=q/p$ and let $X$ be a multiset of $n$
points. Let $\mathcal F$ consist of the convex hulls of all
submultisets of $X$ of size $\lceil\varepsilon n\rceil$, repeating
hulls if needed to reach $p$ members. Any $p$ members have at
least $p\varepsilon n=qn$ incidences with $X$, counted with
multiplicity, so some point of $X$ lies in at least $q$ of them.
Thus $\mathcal F$ has the $(p,q)$-property. Every convex set
containing at least $\varepsilon n$ points of $X$ contains a
member of $\mathcal F$, so any piercing set for $\mathcal F$ is
a weak $\varepsilon$-net for $X$. Taking the worst case over $X$
gives $\HD_d(p,q)\ge f_d(q/p)$ (see Lemma~\ref{lem:average}).}
\[
 \HD_2(p,q)\ge c\,\frac{p}{q}\log\frac{p}{q}.
\]
In particular, $\HD_2(p,q)>p-q+1$ for all $3\le q\le c'\log p$,
provided $c'$ is sufficiently small and $p$ is sufficiently large.

For $d\ge3$, Theorem~\ref{thm:main} gives
$\HD_d(p,q)=p-q+1$ for all $q\ge p^{1-1/\rho_d+\eta}$
and sufficiently large $p$, for every fixed $\eta>0$,
where $\rho_d$ is the exponent in Rubin's weak-net bound
(\ref{eq:rubin}). On the other hand, Bhore, Browne, and Zheng
\cite[Corollary~1]{BBZ26} showed that
\[
 \HD_d(p,q)=\Omega_d\!\left((p/q)^{\gamma_d}\right),
 \qquad
 \gamma_d=\frac12\left\lfloor\frac{\sqrt{8d+25}-1}{2}\right\rfloor .
\]
In particular, for each fixed $d\ge3$, the equality $\HD_d(p,q)=p-q+1$ can hold only
for $q\ge c_d p^{1-1/\gamma_d}$, for a suitable constant $c_d>0$.
We suspect that the actual range is much narrower.

\paragraph{Organization.}
We prove Theorem~\ref{thm:main} in Section~\ref{sec:mainproof}
and Theorem~\ref{thm:correspondence} in Section~\ref{sec:correspondence}.
In Appendix~\ref{sec:nearsharp}, we show how optimizing the parameters
in the proof of \cite{KST18} extends their near-tight upper bound
$p-q+2$ to the wider range stated in~(\ref{eq:nearintro}).

\section{Proof of Theorem~\ref{thm:main}}\label{sec:mainproof}

\subsection{Passing to sequences of convex sets}

The main tool that saves the extra piercing point is passing to sequences of convex sets. A \emph{convex $(\lambda,d)$-body} is an
ordered sequence $(K_1,\ldots,K_\lambda)$ of convex sets in
$\RR^d$; components may be empty. A family of such sequences is
\emph{strongly intersecting} if there is one coordinate $i$ in
which all the components have a common point. It is \emph{weakly
intersecting} if points $z_1,\ldots,z_\lambda$ can be chosen so
that each sequence contains $z_i$ in its $i$th component for at
least one $i$. In the weak condition, the successful coordinate
may depend on the sequence.

For $\lambda\ge1$, put
\[
 h_d(\lambda)=
 \bigl(\lfloor\log_2(\lambda+1)\rfloor+1\bigr)^d.
 \tag{8}\label{eq:h}
\]

\begin{krtheorem}[\cite{KR99}, Theorem 6.1]
If every subfamily of at most $h_d(\lambda)$ convex
$(\lambda,d)$-bodies is strongly intersecting, then the whole
finite family is weakly intersecting.
\end{krtheorem}

For a family with at least $h_d(\lambda)$ members, it suffices to
check subfamilies of exactly that size. Stating the hypothesis
with ``at most'' also handles smaller families without a vacuity
convention. We will use the contrapositive: failure of weak
intersection has a small certificate of failure of strong
intersection.

The point saving follows from the next two lemmas. Say that a
family $\mathcal X$ can be \emph{absorbed} into
$D_1,\ldots,D_\lambda$ if one can choose $z_i\in D_i$ that
together pierce all of $\mathcal X$. The next lemma shows that if $\mathcal X$ cannot be absorbed,
then some `small' subfamily of at most
$h_d(\lambda)=O_d((\log(\lambda+1))^d)$ sets has no common point
in any of the $D_i$.

\begin{lemma}
\label{lem:certificate}
Let $\lambda\ge1$, let $D_1,\ldots,D_\lambda$ be nonempty compact convex sets, and
$\mathcal X$ a family of compact convex sets. If $\mathcal X$
cannot be absorbed into the $D_i$, there is a nonempty
$E\subseteq\mathcal X$ with $|E|\le h_d(\lambda)$ such that
\[
 D_i\cap\bigcap_{A\in E}A=\varnothing
 \qquad\text{for every }i.
 \tag{9}\label{eq:certificate}
\]
\end{lemma}
\begin{proof}
For each $A\in\mathcal X$, form the sequence
\[
 (A\cap D_1,\ldots,A\cap D_\lambda).
\]
A weak transversal of these sequences would give the desired
absorbing points. If a coordinate point lies outside $D_i$, it
pierces no component in that coordinate and can be replaced by
an arbitrary point of $D_i$. Thus weak intersection fails.
Kaiser--Rabinovich gives at most $h_d(\lambda)$ sequences that
fail strong intersection, which is exactly
(\ref{eq:certificate}). If $|\mathcal X|<h_d(\lambda)$, take
$E=\mathcal X$: strong intersection of the entire sequence family
would already imply absorption.
\end{proof}

The following lemma saves the extra piercing point in the proof
of Theorem~\ref{thm:main}, allowing us to improve the bound
$p-q+2$ of \cite{KST18} to the optimal bound $p-q+1$.
\begin{lemma}
\label{lem:absorption}
Suppose $\lambda\ge1$ and a $(\lambda+q,q)$-family is partitioned as
\[
 \F=\{D_1,\ldots,D_\lambda\}\,\dot\cup\,\C\,\dot\cup\,\R,
\]
where the $D_i$ are pairwise disjoint and each of $\C,\R$ is
intersecting. If
\[
 q\ge2h_d(\lambda)+d+1,
 \tag{10}\label{eq:abscondition}
\]
then $\tau(\F)\le\lambda+1$.
\end{lemma}

\begin{proof}
Using one point in each $D_i$ and one for each of $\C$ and $\R$
gives $\lambda+2$ piercing points. If $\C\cup\R$ has a common
point, that point pierces both subfamilies. Otherwise, we show
that the points in the $D_i$ can be chosen to pierce all the
sets in one of the two subfamilies as well. One additional
point then pierces the other subfamily. It remains to consider the case where $\C$ and $\R$ are both
nonempty and $\C\cup\R$ has empty intersection. Suppose for a
contradiction that neither subfamily can be absorbed into the $D_i$.

Lemma~\ref{lem:certificate} supplies disjoint indexed subfamilies
$E_{\C}\subseteq\C$ and $E_{\R}\subseteq\R$, each of size
at most $h_d(\lambda)$, whose respective intersections miss every
$D_i$. Since $\C\cup\R$ is not intersecting, Helly gives
$U\subseteq\C\cup\R$ with $|U|\le d+1$ and empty intersection.
Their union has at most $q$ members by
(\ref{eq:abscondition}). There are at least $q$ members in
$\C\cup\R$, because $|\F|\ge\lambda+q$. Extend the union
to a family $Z\subseteq\C\cup\R$ of exactly $q$ members.

Apply the $(\lambda+q,q)$-property to $\{D_1,\ldots,D_\lambda\}
\cup Z$. An intersecting $q$-tuple cannot use two $D_i$, since
they are disjoint. It cannot use none, since it would be all of
$Z$, including the empty-intersection witness $U$. It therefore
uses exactly one $D_i$ and all but one member of $Z$.

Omitting one member cannot destroy both disjoint certificates
$E_{\C}$ and $E_{\R}$. Hence the $q$-tuple contains a complete
certificate together with $D_i$, contradicting
(\ref{eq:certificate}). This proves the lemma.
\end{proof}

\subsection{The dichotomy}

The dichotomy in \cite{KST18} passes from a missing $(p',q')$-property to a
smaller $(p-p',q-q'+1)$-property outside a witness. We choose a
constant shift that forces the outside family to be intersecting.

\begin{proposition}\label{prop:transfer}
Let $p\ge q+1$, $q\ge2d+1$, and $\lambda=p-q$. If
$q\ge2h_d(\lambda)+d+1$, then
\[
 \HD_d(p,q)\le
 \max\{\HD_d(p-d-1,q-d),\ p-q+1\}.
 \tag{11}\label{eq:transfer}
\]
\end{proposition}
\begin{proof}
Put $b=d+1$ and let $\F$ have the $(p,q)$-property. If it also
has the $(p-b,q-b+1)$-property, the first term in
(\ref{eq:transfer}) applies. Its parameters are admissible:
$q-d\ge d+1$ and their difference is $p-q-1\ge0$.

Otherwise choose $S\subseteq\F$ of size $p-b$ with no
intersecting $(q-b+1)$-tuple. Let $m$ be its largest intersecting
size, so $m\le q-b$, and put $\R=\F\setminus S$.
There are at least $b$ members in $\R$. For any $b$ of them,
denoted by $B$, apply the original $(p,q)$-property to $S\cup B$.
An intersecting $q$-tuple uses at most $q-b$ members of $S$;
therefore it must use all of $B$ and exactly $q-b$ members of $S$.
Consequently $m=q-b$, and every $d+1$ members of $\R$ intersect.
Helly's theorem makes $\R$ intersecting.

Choose intersecting $\C\subseteq S$ of size $q-b$ and put
$\D=S\setminus\C$. Then
\[
 |\D|=(p-b)-(q-b)=p-q=\lambda.
\]
If two members of $\D$ intersect, pierce them together, pierce
the remaining $\lambda-2$ members separately, and add one point
for each of $\C,\R$. This uses $\lambda+1$ points.
Otherwise $\D$ is pairwise disjoint, and
Lemma~\ref{lem:absorption} gives the same bound. This proves
(\ref{eq:transfer}).
\end{proof}

\subsection{The exact bound in the enlarged range}

To complete the proof, we use the following technical update of
the near-sharp estimate of \cite{KST18}, with the weak-net exponent as a parameter.
The proof is given in Appendix~\ref{sec:nearsharp}.

\begin{theorem}\label{thm:nearsharp}
Suppose that for some $s>1$ and every fixed $\sigma>0$,
\[
 f_d(\varepsilon)=O_{d,\sigma}(\varepsilon^{-s-\sigma}).
\]
Then, for every fixed $\eta>0$ and sufficiently large $p$,
\[
 p\ge q\ge d+1,\qquad q\ge p^{1-1/s+\eta}
 \quad\Longrightarrow\quad \HD_d(p,q)\le p-q+2.
 \tag{12}\label{eq:nearsharp}
\]
(In the plane, Rubin's bound gives $p-q+2$ for
$q\ge p^{1/3+\eta}$.)
\end{theorem}

\begin{proof}[Proof of Theorem~\ref{thm:main}]
The case $p=q$ follows from Helly, so assume $p\ge q+1$.
Write $\beta=1-1/\rho_d$, and apply
Theorem~\ref{thm:nearsharp} to
\[
 P=p-d-1,\qquad Q=q-d,
\]
with exponent margin $\eta/2$. For sufficiently large $p$,
these parameters satisfy $P\ge Q\ge d+1$ and
\[
 Q\ge p^{\beta+\eta}-d
 \ge\tfrac12p^{\beta+\eta}
 \ge P^{\beta+\eta/2}.
\]
Thus
\[
 \HD_d(P,Q)\le P-Q+2=p-q+1.
\]
Also $2h_d(p-q)+d+1=O_d(\log^d p)=o(q)$, so
Proposition~\ref{prop:transfer} applies and gives the required
upper bound. The lower bound is (\ref{eq:standardlower}).
\end{proof}

\section{Proof of Theorem~\ref{thm:correspondence}}
\label{sec:correspondence}

The connection is easiest to see in reverse. If a family consists
of convex sets that each contain an $\varepsilon$ fraction of a point
set, any $p$ members have at least $\varepsilon p$ incidences on
average at a point. Thus a sufficiently strong exact $(p,q)$
theorem immediately produces a weak net. The other direction
uses a net to control the error in the dense branch of the
dichotomy, while the finite transfer removes the last extra point.

For each integer $p\ge d+1$, let $Q_d(p)$ be the least
$q_0\in\{d+1,\ldots,p\}$ such that $\HD_d(p,q)=p-q+1$
for every $q_0\le q\le p$. For $t\ge1$, define
\[
 u(t)=\inf\{\varepsilon\in(0,1]:f_d(\varepsilon)\le t\},
\]
and set $u(t)=1$ for $0<t<1$.

The equivalence in Theorem~\ref{thm:correspondence} follows
from the following quantitative estimates. For some $C_d>0$
and all sufficiently large $p$,
\[
 p u(p)\le Q_d(p)\le C_d p u(p)\log^2(p+e).
 \tag{13}\label{eq:comparison}
\]
For every $0<\varepsilon\le1$,
\[
 f_d(\varepsilon)\le
 \min\{p\in\mathbb N:p\ge d+1,\ \varepsilon p\ge Q_d(p)\}.
 \tag{14}\label{eq:inverseQ}
\]
The function $u$ is nonincreasing. At a jump point the infimum
need not be attained, but monotonicity always gives
\[
 f_d(\min\{1,2u(t)\})\le t\qquad(t\ge1).
 \tag{15}\label{eq:inverseproperty}
\]

\subsection{From an exact range to weak nets}

Before getting into the following lemma, we note that finite point sets and finite point multisets have the same
optimal weak $\varepsilon$-net bound. Indeed, perturb the
occurrences of a multiset $X$ to distinct points and pass to a
convergent subsequence of the corresponding nets in a fixed
closed ball as the perturbations tend to zero. The limiting
net has the same size bound and meets the convex hull of every
$\lceil\varepsilon |X|\rceil$-element submultiset of $X$.
We recall the averaging observation used in \cite{KST18}, just
before Proposition 2.8.

\begin{lemma}[Averaging]\label{lem:average}
If $p\ge q\ge d+1$ and $\varepsilon p\ge q$, then
$f_d(\varepsilon)\le\HD_d(p,q)$. In particular,
$\HD_d(p,q)\le p$ implies $f_d(\varepsilon)\le p$.
\end{lemma}
\begin{proof}
Given an $n$-point set $X$, take the finite family of convex hulls
of all its subsets of size $\lceil\varepsilon n\rceil$. Piercing
these hulls pierces every $\varepsilon$-heavy convex set. If there
are fewer than $p$ hulls, repeat members until the indexed
family has at least $p$; this does not change its piercing number.
Any $p$ members have at least $p\varepsilon n\ge qn$ incidences
with $X$, so some point belongs to at least $q$ of them.
The family has the $(p,q)$-property and can be pierced by at
most $\HD_d(p,q)$ points.
\end{proof}

Taking $q=Q_d(p)$ proves (\ref{eq:inverseQ}). Taking
$\varepsilon=Q_d(p)/p$ gives $f_d(Q_d(p)/p)\le p$, and hence
\[
 p u(p)\le Q_d(p).
 \tag{16}\label{eq:necessary}
\]
More explicitly, if a proved exact range has threshold
$Q_d(p)\le A p^a\log^b(p+e)$, with $0\le a<1$ and $b\ge0$,
then
\[
 f_d(\varepsilon)=O\!\left(
 \varepsilon^{-1/(1-a)}\log^{b/(1-a)}(2/\varepsilon)\right).
 \tag{17}\label{eq:explicitreverse}
\]
Indeed, choose $p$ to be a sufficiently large constant times
$\varepsilon^{-1/(1-a)}\log^{b/(1-a)}(2/\varepsilon)$, rounded up.
Then $Q_d(p)\le\varepsilon p$, and apply Lemma~\ref{lem:average}.

\subsection{A sufficient threshold for a general net function}

We first record a scaling fact that will let us simplify the final
threshold even when $f_d$ has irregular growth. By the observation preceding Lemma~\ref{lem:average}, we may use a
finite point set to witness $f_d(\varepsilon)$ in the proof
of the next lemma.

\begin{lemma}[Replication]\label{lem:replication}
For every positive integer $k$,
\[
 f_d(\varepsilon/k)\ge kf_d(\varepsilon).
\]
Consequently, for $1\le x\le y$,
\[
 xu(x)\le2yu(y).
 \tag{18}\label{eq:inverseprofile}
\]
\end{lemma}
\begin{proof}
Take a point set witnessing $f_d(\varepsilon)$, together with its
finite family of heavy convex hulls. Place $k$ translated and
scaled copies in disjoint balls. Each local hull is
$(\varepsilon/k)$-heavy in the union, and a piercing point cannot
serve two balls. At least $kf_d(\varepsilon)$ points are needed.
A witnessing finite set exists because the optimal bound is
integer-valued and finite.

For the second assertion, the case $u(y)=1$ is immediate.
Otherwise put $k=\lceil y/x\rceil$ and choose
$u(y)<\varepsilon'\le1$ with $f_d(\varepsilon')\le y$.
If $k\varepsilon'\le1$, the first assertion gives
$f_d(k\varepsilon')\le y/k\le x$, so $u(x)\le k\varepsilon'$.
If $k\varepsilon'>1$, the same inequality follows from $u(x)\le1$.
Letting $\varepsilon'\downarrow u(y)$ yields
$xu(x)\le x\lceil y/x\rceil u(y)\le2yu(y)$.
\end{proof}

We use the following form of KST Proposition 2.8:
\[
 \HD_d(P,V)\le P-V+f_d(\varepsilon_0)+B_d,
 \qquad
 \varepsilon_0\ge\frac{c_*V}{P\log(2P/V)},
 \tag{19}\label{eq:KSTinput}
\]
when $P\ge V\ge d+1$ and $V\ge\log P$. Here
$c_*,B_d>0$ depend only on $d$, and $\varepsilon_0\in(0,1]$.
The regularized logarithm handles $P/V$ close to one. The bounded
term accommodates the final small-subfamily bookkeeping and will
have no effect below.

\begin{proposition}[An integral threshold]\label{prop:integral}
There are constants $\xi,K>0$ such that, with
\[
 H(x)=K\left[
 1+\log^d(x+e)+\int_1^x u(\xi t)\log(t+e)\,dt
 \right],\qquad x\ge1,
 \tag{20}\label{eq:H}
\]
the exact bound $\HD_d(p,q)=p-q+1$ holds whenever
$p\ge q\ge d+1$ and $q\ge H(p)$.
\end{proposition}
\begin{proof}
Let $\alpha=(d-1)/d$, choose $\theta\in(\alpha,1)$, and put
$\kappa=\theta-\alpha>0$. Fix $0<\xi<\kappa/8$.
We will choose $K$ sufficiently large and use strong induction
on $p$. Since $f_d(\varepsilon)\ge\lfloor1/\varepsilon\rfloor
\ge1/(2\varepsilon)$, we have
\[
 u(t)\ge1/(2t)\quad(t\ge1).
 \tag{21}\label{eq:inversebasic}
\]
Also $u(t)\to0$, since $f_d(\varepsilon)$ is finite for every
$\varepsilon>0$.

For the induction step, first consider a pair
\[
 N\le P<p,\quad P\ge Q\ge d+1,\quad Q\ge H(P)-d,
 \tag{22}\label{eq:auxpair}
\]
where $N$ will be fixed sufficiently large. If
$Q>\alpha P+1$, the classical theorem gives the exact bound.
Otherwise put
\[
 t=\lfloor\theta P\rfloor,\qquad \ell=\lceil H(t)\rceil,
 \qquad v=Q-\ell.
\]
Take $N$ large enough that $t\ge d+1$ and $P-t\ge1$.
The integral in (\ref{eq:H}) and monotonicity of $u$ give
\[
 H(P)-H(t)
 \ge K(1-\theta)P\,u(\xi P)\log(\theta P+e).
\]
Consequently, for some $c_\theta>0$ and sufficiently large $N$,
\[
 v\ge H(P)-H(t)-d-1
 \ge Kc_\theta P\,u(\xi P)\log(P+e).
 \tag{23}\label{eq:gap}
\]
The right side is at least a constant times $K\log P$, by
(\ref{eq:inversebasic}). After increasing $K$, it follows that
$v\ge\max\{d+1,\log P\}$. Since $v>0$,
$\ell<Q\le\alpha P+1$, and therefore
\[
 t-\ell\ge\kappa P-2\ge\kappa P/2.
 \tag{24}\label{eq:linearsaving}
\]

Let $\F$ have the $(P,Q)$-property. If it also has the
$(P-t,v)$-property, necessarily $v\le P-t$, so
(\ref{eq:KSTinput}) applies. Its weak-net parameter satisfies
\[
 \varepsilon_0\ge
 \frac{c_*v}{(P-t)\log(2(P-t)/v)}
 \ge c'_d K u(\xi P).
\]
Choose $K$ so that $c'_dK\ge2$, and take $N$ large enough
that $2u(\xi P)\le1$. By (\ref{eq:inverseproperty}),
$f_d(\varepsilon_0)\le\xi P$. Hence
\[
 \tau(\F)\le P-Q-(t-\ell)+\xi P+B_d\le P-Q+1
\]
using (\ref{eq:linearsaving}) and a sufficiently large $N$.

Otherwise choose $S\subseteq\F$ of size $P-t$ with no
intersecting $v$-tuple. If $v>P-t$, any such $S$ works.
Let $m$ be its largest intersecting size. Then $m\le v-1$,
and $\R=\F\setminus S$ has the $(t,Q-m)$-property:
apply the original property to $S$ and any $t$ members of
$\R$. There are at least $t$ such members, and this argument
also forces $Q-m\le t$. Moreover
\[
 Q-m\ge\ell+1\ge H(t).
\]
The strong induction hypothesis at $t<p$ gives
$\tau(\R)\le t-Q+m+1$. Piercing a maximum intersecting
subfamily of $S$ once and all its other members individually
costs at most $P-t-m+1$. Thus in this case
\[
 \tau(\F)\le P-Q+2.
 \tag{25}\label{eq:auxplus2}
\]

Now let $q\ge H(p)$. If $q=p$, Helly suffices. Otherwise,
the term $K\log^d(p+e)$ ensures all conditions of
Proposition~\ref{prop:transfer}, for $K$ large enough.
The shifted pair
\[
 P=p-d-1,\qquad Q=q-d\ge H(p)-d\ge H(P)-d
\]
satisfies (\ref{eq:auxpair}), and its bound
$P-Q+2=p-q+1$ proves the induction step.

All lower thresholds for $P$ above can be chosen uniformly for
increasing $K$: (\ref{eq:inversebasic}) makes the expression
in (\ref{eq:gap}) grow at least as $K\log P$, and every other
increase in $K$ strengthens the estimates. Fix such an $N$ and
finally enlarge $K$ so that $H(p)>p$ for every integer
$p\le N+d+1$. The induction base is then vacuous. The lower
bound (\ref{eq:standardlower}) gives equality.
\end{proof}

\subsection{Completing the comparison and its consequences}

For $\xi t\ge1$ and $t\le p$, Lemma~\ref{lem:replication}
gives $u(\xi t)\le2p u(p)/(\xi t)$. The remaining bounded
interval contributes only a constant. Therefore
\[
 \int_1^p u(\xi t)\log(t+e)\,dt
 \le C_d p u(p)\log^2(p+e).
 \tag{26}\label{eq:intbound}
\]
By (\ref{eq:BMN}),
\[
 p u(p)\ge c_d\log^{d-1}(p+e)
 \tag{27}\label{eq:BMNinverse}
\]
for large $p$. Indeed, if
$\varepsilon<c\log^{d-1}p/p$, then $\log(1/\varepsilon)\ge\frac12\log p$
and (\ref{eq:BMN}) exceeds $p$ when $c$ is small enough.
In particular the term $\log^d(p+e)$ in (\ref{eq:H}) is
absorbed by the right side of (\ref{eq:intbound}).
Proposition~\ref{prop:integral} now proves the upper bound in
(\ref{eq:comparison}); (\ref{eq:necessary}) proves the lower
bound, and Lemma~\ref{lem:average} proves (\ref{eq:inverseQ}).

We verify the exponent equivalence in
Theorem~\ref{thm:correspondence}. If
$f_d(\varepsilon)=O_\sigma(\varepsilon^{-s-\sigma})$, then
$u(p)=O_\sigma(p^{-1/(s+\sigma)})$ and
\[
 Q_d(p)=O_\sigma\!\left(
 p^{1-1/(s+\sigma)}\log^2(p+e)\right).
\]
For any desired $\eta>0$, choose $\sigma>0$ small
enough and absorb the logarithm into the remaining power margin.
Conversely, if the exact range holds above
$p^{1-1/s+\eta}$ for every positive margin, choose
$a>1-1/s$ so close to it that $1/(1-a)<s+\sigma$.
Formula (\ref{eq:explicitreverse}) with $b=0$ gives
$f_d(\varepsilon)=O(\varepsilon^{-1/(1-a)})$. This proves the converse,
including $s=1$ by choosing a sufficiently small $a>0$.

Thus if the optimal polynomial exponents are denoted by $s_*$
for weak nets and $a_*$ for the sharp range, with arbitrary
positive power margins allowed, then
\[
 a_*=1-1/s_*,\qquad s_*=1/(1-a_*).
 \tag{28}\label{eq:exponents}
\]
For $d\ge3$, the lower bound of Bhore, Browne, and Zheng
\cite[Theorem~1]{BBZ26} gives
$f_d(\varepsilon)=\Omega_d(\varepsilon^{-\gamma_d})$, and hence
$u(p)=\Omega_d(p^{-1/\gamma_d})$. By (\ref{eq:necessary}), it
follows that $Q_d(p)=\Omega_d(p^{1-1/\gamma_d})$.
Together with Rubin's upper bound (\ref{eq:rubin}) and the
exponent equivalence, this gives
\[
 1-\frac{1}{\gamma_d}\le a_*\le 1-\frac{1}{\rho_d}
 \qquad(d\ge3).
\]

If an asymptotic power exists, the pointwise statement is
$f_d(\varepsilon)=\varepsilon^{-s+o(1)}$ if and only if
$Q_d(p)=p^{1-1/s+o(1)}$, by monotone inversion and
(\ref{eq:comparison}). This does not identify the logarithmic
factors.

For example, if a valid upper bound is
$f_d(\varepsilon)\le A\varepsilon^{-s}\log^k(2/\varepsilon)$, with $s>1$,
then evaluating the integral in (\ref{eq:H}) yields
\[
 Q_d(p)=O\!\left(
 p^{1-1/s}\log^{1+k/s}(p+e)\right).
\]
In dimension $2$, if
$f_2(\varepsilon)=O(\varepsilon^{-1}\log^k(2/\varepsilon))$
for some fixed $k\ge0$, then evaluating the same integral gives
$Q_2(p)=O(\log^{k+2}(p+e))$.

\appendix
\section{A technical update of the near-sharp bound}
\label{sec:nearsharp}

We give the details of the technical update to the proof of the
$p-q+2$ bound in \cite[Section~2.3]{KST18}. In that proof, the
dimension $d$ enters both as the exponent in the weak-net bound
and in the classical Hadwiger--Debrunner condition. We keep
the geometric split and use a general weak-net exponent $s$ in
the bootstrapping. Taking $s=\rho_d$ using Rubin's bounds
\cite{RubinPlane,RubinHigher} gives (\ref{eq:nearintro}).

Assuming $f_d(\varepsilon)=O_\sigma(\varepsilon^{-s-\sigma})$ for
every $\sigma>0$, our goal is to prove
\[
 \HD_d(p,q)\le p-q+2
 \quad\text{whenever }q\ge p^{1-1/s+\eta},
\]
for every fixed $\eta>0$ and sufficiently large $p$.

Fix $r>s$. From (\ref{eq:KSTinput}) there is a constant $C$ such
that
\[
 \HD_d(P,Q)\le P-Q+C(P/Q)^r
 \qquad(Q\ge\log P).
 \tag{29}\label{eq:polynomialinput}
\]
To see this, insert a weak-net exponent strictly between $s$
and $r$, then absorb the logarithmic factor into the remaining
positive power. Constants may depend on $d,s,r$.

\subsection{The geometric split from the KST proof}

Let $\alpha=(d-1)/d$. The case $q>\alpha p+1$ is classical,
so suppose $q\le\alpha p+1$. For large $q$, put
\[
 h=\left\lfloor\frac q{d-1}\right\rfloor,
 \qquad j=\left\lceil\frac qd\right\rceil,
 \qquad p'=p-h,\quad q'=(d-1)j-d.
 \tag{30}\label{eq:geosplit}
\]
These parameters are admissible and satisfy $p'\asymp_d p$,
$q'\asymp_d q$, as well as
\[
 p'-q'\le p-q-\frac{q}{d(d-1)}+d+1.
 \tag{31}\label{eq:geosaving}
\]
Indeed, $h\ge q/(d-1)-1$ and $j\ge q/d$ give the displayed
inequality. Also $q\le\alpha p+1$ gives
$h\le p/d+O_d(1)$, while
$p'-q'\ge p-(1+1/[d(d-1)])q+O_d(1)$.
In the nonclassical range this last expression is at least
$((d-1)/d^2)p+O_d(1)$, so $p'\ge q'$ for large $p$.

\begin{lemma}[The sparse alternative in the geometric split]
\label{lem:geometric}
If a $(p,q)$-family does not have the $(p',q')$-property in
(\ref{eq:geosplit}), it is pierced by at most $p-q+2$ points.
\end{lemma}
\begin{proof}
Choose $S$ of size $p-h$ with no intersecting $q'$-tuple.
Write its maximum intersecting size as
$m=(d-1)j-t$, where $t\ge d+1$.
The remaining family has the $(h,q-m)$-property, by applying
the original property to $S$ together with any $h$ outside
members. This forces $q-m\le h$. Also $j\le q/d+1$ and
$t\ge d+1$ give $q-m\ge q/d+2$, so $q-m\ge d+1$ for large $q$.

The parameters meet the classical condition, since
\[
 \begin{aligned}
 (d-1)h
 &\le q
 <q+d(t-d)\\
 &\le d\bigl(q-(d-1)j+t-1\bigr)=d(q-m-1).
 \end{aligned}
\]
Here the last inequality follows from $j\le q/d+1$.
Therefore the remaining family is pierced by at most
$h-q+m+1$ points. Piercing a maximum intersecting subfamily
of $S$ once and all its other members separately uses at most
$p-h-m+1$ points. The sum is $p-q+2$.
\end{proof}

\subsection{Bootstrapping a bounded excess}

Put
\[
 a_k=\frac{r^k(r-1)}{r^{k+1}-1},\qquad k\ge1.
 \tag{32}\label{eq:ak}
\]
These exponents satisfy
$a_1=r/(r+1)$,
$a_{k+1}=ra_k/(1+ra_k)$, and
$a_k\downarrow1-1/r$.

\begin{proposition}\label{prop:boundedexcess}
For every fixed integer $k\ge1$ and $\eta>0$, and all
sufficiently large $p$,
\[
 q\ge p^{a_k+\eta}
 \quad\Longrightarrow\quad \HD_d(p,q)\le p-q+k+1.
 \tag{33}\label{eq:boundedexcess}
\]
\end{proposition}
\begin{proof}
For $k=1$, use (\ref{eq:geosplit}). Its sparse alternative
is handled by Lemma~\ref{lem:geometric}. In the dense
alternative, (\ref{eq:polynomialinput}) and
(\ref{eq:geosaving}) give
\[
 \tau(\F)\le p-q-\frac{q}{d(d-1)}+d+1+O_d((p/q)^r).
\]
The error is $o(q)$ when $q\ge p^{r/(r+1)+\eta}$, since
$(p/q)^r/q\le p^{-(r+1)\eta}$. This proves the initial case.

Suppose the result is known for $k$. It is enough to prove the
next case for a fixed exponent $b$ with
$a_{k+1}<b<a_k$. Once the assertion holds for such a $b$, it
also holds for any larger threshold exponent, by inclusion of
the corresponding ranges. Choose $A>a_k$, sufficiently close to $a_k$, so that
\[
 A<1,\qquad \frac{rA}{1+rA}<b.
\]
If $q\ge p^A$, the induction hypothesis already gives the
claim. Otherwise assume $p^b\le q<p^A$ and put
\[
 t=\left\lfloor(q/4)^{1/A}\right\rfloor,
 \qquad w=\lfloor q/2\rfloor.
\]
Then $t<p/4$, $q=o(t)$, and $t^A\le q/4$. Moreover
\[
 \frac{(p/q)^r}{t}
 =O\!\left(\frac{p^r}{q^{r+1/A}}\right)
 \le O\!\left(p^{r-b(r+1/A)}\right)=o(1).
 \tag{34}\label{eq:bootstrapdominance}
\]

If $\F$ has the $(p-t,q-w)$-property, these parameters are
admissible for large $p$, with $q-w\gg\log p$.
Formula (\ref{eq:polynomialinput}) yields
\[
 \tau(\F)\le p-q-t+w+O((p/q)^r)<p-q.
\]
Otherwise choose $S$ of size $p-t$ with no intersecting
$(q-w)$-tuple. If $m$ is its maximum intersecting size, then
$m\le q-w-1$, and $\F\setminus S$ has the $(t,q-m)$-property.
The original property forces $q-m\le t$, while
\[
 q-m\ge w+1>q/2\ge t^A.
\]
Since $A>a_k$, the induction hypothesis applies and gives
\[
 \tau(\F\setminus S)\le t-q+m+k+1.
\]
Adding $p-t-m+1$ points for $S$ gives $p-q+k+2$.
This proves the induction step. Every parameter grows as a
positive power, so all fixed lower thresholds and the integer
roundings above are harmless once $p$ is sufficiently large.
\end{proof}

\subsection{Returning to an excess of two}

\begin{proof}[Proof of Theorem~\ref{thm:nearsharp}]
Fix $\eta>0$ and choose $r>s$ sufficiently close to
$s$ that
\[
 1-1/r<\beta:=1-1/s+\eta.
\]
We may suppose $\beta<1$; if $\beta=1$, Helly handles the only
possible pair $q=p$, and if $\beta>1$ the range is empty for
$p>1$. Choose a fixed $k$ so that $a_k<\beta$.
Suppose $q\ge p^\beta$. The classical case is already done,
so use (\ref{eq:geosplit}) once more.

Its sparse alternative gives $p-q+2$ by
Lemma~\ref{lem:geometric}. In the dense alternative,
$q'\asymp_d q$ and $p'\asymp_d p$, so
Proposition~\ref{prop:boundedexcess} applies to $(p',q')$
with a fixed positive exponent margin. Therefore
\[
 \begin{aligned}
 \tau(\F)&\le p'-q'+k+1\\
 &\le p-q-\frac{q}{d(d-1)}+d+k+2\\
 &\le p-q+2
 \end{aligned}
\]
for sufficiently large $p$. This proves the theorem using only
the KST dichotomy, the weak-net input, and the classical
Hadwiger--Debrunner theorem.
\end{proof}

\paragraph{Use of AI.}
The ideas in this paper were developed through discussions between
the author and OpenAI's GPT-6.1 Sol. The initial draft was written
using ChatGPT and subsequently revised by the author.


\begin{thebibliography}{99}
\bibitem{HD57}
H.~Hadwiger and H.~Debrunner,
\"Uber eine Variante zum Hellyschen Satz,
\emph{Arch. Math.} \textbf{8} (1957), 309--313.
\href{https://doi.org/10.1007/BF01898794}{doi:10.1007/BF01898794}.

\bibitem{AK92}
N.~Alon and D.~J.~Kleitman,
Piercing convex sets and the Hadwiger--Debrunner $(p,q)$-problem,
\emph{Adv. Math.} \textbf{96} (1992), 103--112.
\href{https://doi.org/10.1016/0001-8708(92)90052-M}{doi:10.1016/0001-8708(92)90052-M}.

\bibitem{AK97}
N.~Alon and D.~J.~Kleitman,
A purely combinatorial proof of the Hadwiger--Debrunner $(p,q)$
conjecture, \emph{Electron. J. Combin.} \textbf{4} (1997), no.~2, R1.
\href{https://doi.org/10.37236/1316}{doi:10.37236/1316}.

\bibitem{ABFK92}
N.~Alon, I.~B\'ar\'any, Z.~F\"uredi, and D.~J.~Kleitman,
Point selections and weak $\varepsilon$-nets for convex hulls,
\emph{Combin. Probab. Comput.} \textbf{1} (1992), 189--200.
\href{https://doi.org/10.1017/S0963548300000225}{doi:10.1017/S0963548300000225}.

\bibitem{CEGG+95}
B.~Chazelle, H.~Edelsbrunner, M.~Grigni, L.~J.~Guibas,
M.~Sharir, and E.~Welzl,
Improved bounds on weak $\varepsilon$-nets for convex sets,
\emph{Discrete Comput. Geom.} \textbf{13} (1995), 1--15.
\href{https://doi.org/10.1007/BF02574025}{doi:10.1007/BF02574025}.

\bibitem{KST18}
C.~Keller, S.~Smorodinsky, and G.~Tardos,
Improved bounds on the Hadwiger--Debrunner numbers,
\emph{Israel J. Math.} \textbf{225} (2018), 925--945.
\href{https://arxiv.org/abs/1512.04026v3}{arXiv:1512.04026v3}.

\bibitem{KS20}
C.~Keller and S.~Smorodinsky,
From a $(p,2)$-theorem to a tight $(p,q)$-theorem,
\emph{Discrete Comput. Geom.} \textbf{63} (2020), 821--847.
\href{https://doi.org/10.1007/s00454-018-0048-3}{doi:10.1007/s00454-018-0048-3}.
\href{https://arxiv.org/abs/1712.04552}{arXiv:1712.04552}.


\bibitem{KR99}
T.~Kaiser and Y.~Rabinovich,
Intersection properties of families of convex $(n,d)$-bodies,
\emph{Discrete Comput. Geom.} \textbf{21} (1999), 275--287.
\href{https://doi.org/10.1007/PL00009421}{doi:10.1007/PL00009421}.

\bibitem{RubinPlane}
N.~Rubin,
An improved bound for weak epsilon-nets in the plane,
\emph{J. ACM} \textbf{69} (2022), no.~5, Article 32, 1--35.
\href{https://doi.org/10.1145/3555985}{doi:10.1145/3555985}.
\href{https://arxiv.org/abs/1808.02686v2}{arXiv:1808.02686v2}.

\bibitem{RubinHigher}
N.~Rubin,
Stronger bounds for weak epsilon-nets in higher dimensions,
\emph{Discrete Comput. Geom.} \textbf{76} (2026), 2612--2699.
\href{https://doi.org/10.1007/s00454-026-00878-y}{doi:10.1007/s00454-026-00878-y}.

\bibitem{BMN11}
B.~Bukh, J.~Matou\v{s}ek, and G.~Nivasch,
Lower bounds for weak epsilon-nets and stair-convexity,
\emph{Israel J. Math.} \textbf{182} (2011), 199--228.
\href{https://arxiv.org/abs/0812.5039v2}{arXiv:0812.5039v2}.

\bibitem{KS21}
C.~Keller and S.~Smorodinsky,
A new lower bound on Hadwiger--Debrunner numbers in the plane,
\emph{Israel J. Math.} \textbf{244} (2021), 649--680.
\href{https://doi.org/10.1007/s11856-021-2185-2}{doi:10.1007/s11856-021-2185-2}.

\bibitem{BBZ26}
S.~Bhore, R.~Browne, and D.~W.~Zheng,
Strong lower bounds for weak epsilon-nets,
\href{https://arxiv.org/abs/2610.40982v3}{arXiv:2610.40982v3},
2026.
\end{thebibliography}
\end{document}